\documentclass[draft]{amsart} \usepackage{amssymb}
\newtheorem{thm}{Theorem}[section] \newtheorem{cor}[thm]{Corollary}
\newtheorem{conj}[thm]{Conjecture} \newtheorem{lem}[thm]{Lemma}
\newtheorem{prop}[thm]{Proposition} \theoremstyle{remark}
\newtheorem{rem}[thm]{Remark} \newtheorem{ex}[thm]{Example}
\newtheorem*{question}{Question} \theoremstyle{definition}
\newtheorem{defn}[thm]{Definition} \numberwithin{equation}{section}

\renewcommand{\frak}{\mathfrak}

\newcommand{\bbold}{\mathbb} \newcommand{\Higman}{{\operatorname{H}}}
\newcommand{\lex}{{\operatorname{lex}}}
\newcommand{\Fin}{\operatorname{Fin}} \newcommand{\cal}{\mathcal}
 
\newcommand{\lm}{\operatorname{lm}}
\newcommand{\lc}{\operatorname{lc}}
\newcommand{\lt}{\operatorname{lt}}
\newcommand{\sort}{\operatorname{sort}}

\renewcommand\hat{\widehat}

\newcommand{\<}{\langle} \renewcommand{\>}{\rangle}

\def \N { {\bbold N} }  \def \Q { {\bbold Q} }

\def \i{{\boldsymbol{i}}}

\def \u{{\boldsymbol{u}}}
\def \v{{\boldsymbol{v}}}

\title[Finite generation of symmetric ideals]{Finite generation of
  symmetric ideals}

\author{Matthias Aschenbrenner} \address{Department of Mathematics,
  Statistics, and Computer Science,
University of Illinois at Chicago, Chicago, IL 60607.}%
\email{maschenb@math.uic.edu}%

\author{Christopher J. Hillar} 
\address{Department of Mathematics, Texas A\&M University, College Station, TX 77843}
\email{chillar@math.tamu.edu}

\thanks{The first author was partially supported by the National
  Science Foundation Grant DMS 03-03618.  The work of the second
  author was supported under a National Science Foundation Graduate
  Research Fellowship.}


\dedicatory{In memoriam Karin Gatermann
  \textup{(}1965--2005\textup{)}.}

\date{April 2005}

\begin{document}

\begin{abstract}
  Let $A$ be a commutative Noetherian ring, and let $R = A[X]$ be the
  polynomial ring in an infinite collection $X$ of indeterminates over
  $A$.  Let ${\mathfrak S}_{X}$ be the group of permutations of $X$.  The
  group ${\mathfrak S}_{X}$ acts on $R$ in a natural way, and this in
  turn gives $R$ the structure of a left module over the left group ring
  $R[{\mathfrak S}_{X}]$.  We prove that all ideals of $R$ invariant
  under the action of ${\mathfrak S}_{X}$ are finitely generated as
  $R[{\mathfrak S}_{X}]$-modules.  The proof involves introducing a
  certain well-quasi-ordering on monomials and developing
  a theory of Gr\"obner bases and reduction in this setting.  
  We also consider the concept of an invariant
  chain of ideals for finite-dimensional polynomial rings and relate
  it to the finite generation result mentioned above.  Finally, a
  motivating question from chemistry is presented, with the above
  framework providing a suitable context in which to study it.
\end{abstract}

\maketitle

\section{Introduction}
A pervasive theme in invariant theory is that of finite generation.  A
fundamental example is a theorem of Hilbert stating that the invariant
subrings of finite-dimensional polynomial algebras over finite groups are
finitely generated \cite[Corollary 1.5]{eisenbud}.  In this article,
we study invariant ideals of infinite-dimensional polynomial rings.
Of course, when the number of indeterminates is finite, Hilbert's
basis theorem tells us that any ideal (invariant or not) is finitely
generated.  

Our setup is as follows.  Let $X$ be an infinite collection of
indeterminates, and let ${\mathfrak S}_{X}$ be the group of
permutations of $X$.  Fix a commutative Noetherian ring $A$ and let $R
= A[X]$ be the polynomial ring in the indeterminates $X$.  The group
${\mathfrak S}_{X}$ acts naturally on $R$: if $\sigma \in {\mathfrak
  S}_{X}$ and $f\in A[x_1,\dots,x_n]$, where $x_i\in X$, then
\begin{equation}\label{simpleaction}
\sigma f(x_1,x_2,\ldots,x_n) = f(\sigma x_1,\sigma x_2,\dots,\sigma x_n)\in R.
\end{equation}
Let $R[{\mathfrak S}_{X}]$ be the left group ring associated to ${\mathfrak S}_{X}$ and $R$.  
This ring is the set of all finite linear combinations, 
\[R[{\mathfrak S}_{X}] = \left\{ \sum_{i=1}^m r_i \sigma_i : r_i \in R, \sigma_i \in {\mathfrak S}_{X}\right\}.\]
Multiplication is given by $f\sigma\cdot g\tau = fg(\sigma\tau)$ for $f,g\in R$,
$\sigma,\tau\in {\mathfrak S}_{X}$, and extended by linearity.
The action (\ref{simpleaction}) allows us
to endow $R$ with the structure of a left $R[{\mathfrak S}_{X}]$-module in the natural way.

An ideal $I \subseteq R$ is called \textit{invariant under ${\mathfrak S}_{X}$}
(or simply \textit{invariant}) if \[ {\mathfrak S}_{X}I := \{\sigma f
: \sigma \in {\mathfrak S}_{X}, f \in I\} \subseteq I.\] Notice that
invariant ideals are simply the $R[{\mathfrak S}_{X}]$-submodules of
$R$.  We may now state our main result.

\begin{thm}\label{onevarfinitegenthm}
   Every ideal of $R = A[X]$ invariant under ${\mathfrak
    S}_{X}$ is finitely generated as an $R[ {\mathfrak
    S}_{X}]$-module.  \textup{(}Stated more succinctly, $R$ is a
  Noetherian $R[{\mathfrak S}_{X}]$-module.\textup{)}
\end{thm}

This result is motivated by finiteness questions in
chemistry \cite{Ruch1, Ruch2, Ruch3} and algebraic statistics \cite{SturmSull}
involving chains of invariant ideals $I_k$ ($k = 1,2,\ldots$)
inside finite-dimensional polynomial rings $R_k$.  Section \ref{chemmotivation}
contains a discussion.

For the purposes of this work, we will use the following notation.
Let $B$ be a ring and let $G$ be a subset of a $B$-module $M$.  Then
$\<f: f \in G \>_{B}$ will denote the $B$-submodule of $M$ generated
by elements of $G$.

\begin{ex}
  Suppose that $X=\{x_1,x_2,\dots\}$.  The invariant ideal $I =
  \<x_1,x_2,\ldots \>_{R}$ is clearly not finitely generated over $R$;
  however, it does have the compact representation $I = \<x_1
  \>_{R[{\mathfrak S}_{X}]}.$
\end{ex}

The outline of this paper is as follows.  In Section \ref{partorder},
we define a partial order on monomials and show that it can be used to
obtain a well-quasi-ordering of the monomials in $R$.  Section
\ref{mainthmproof} then goes on to detail our proof of Theorem
\ref{onevarfinitegenthm}, using the main result of Section
\ref{partorder} in a fundamental way.  In the penultimate section, we
discuss a relationship between invariant ideals of $R$ and chains of
increasing ideals in finite-dimensional polynomial rings.  The notions
introduced there provide a suitable framework for studying a problem
arising from chemistry, the subject of the final section of this
article.

\section{The Symmetric Cancellation Ordering}\label{partorder}

We begin this section by briefly recalling some basic order-theoretic
notions.  We also discuss some fundamental results due to Higman and
Nash-Williams and some of their consequences. We define the ordering
mentioned in the section heading and give a sufficient condition for
it to be a well-quasi-ordering; this is needed in the proof of
Theorem~\ref{onevarfinitegenthm}.

\subsection{Preliminaries}
A {\it quasi-ordering} on a set $S$ is a binary relation $\leq$ on $S$
which is reflexive and transitive.  A {\it quasi-ordered set} is a
pair $(S,\leq)$ consisting of a set $S$ and a quasi-ordering $\leq$ on
$S$. When there is no confusion, we will omit $\leq$ from the
notation and simply call $S$ a quasi-ordered set.  If in addition the
relation $\leq$ is \textit{anti-symmetric} ($s \leq t \ \wedge \ t
\leq s \Rightarrow s = t$, for all $s,t\in S$), then $\leq$ is called
an {\it ordering} (sometimes also called a {\it partial ordering}\/) on the
set $S$.  The {\it trivial} ordering on $S$ is given by $s\leq
t\Longleftrightarrow s=t$ for all $s,t\in S$.  A quasi-ordering $\leq$
on a set $S$ induces an ordering on the set $S/{\sim}=\{s/{\sim}:s\in
S\}$ of equivalence classes of the equivalence relation $s\sim t
\Longleftrightarrow s\leq t \ \wedge\ t\leq s$ on $S$.  If $s$ and $t$
are elements of a quasi-ordered set, we write as usual $s\leq t$ also
as $t\geq s$, and we write $s<t$ if $s\leq t$ and $t\not\leq s$.

A map $\varphi\colon S\to T$ between quasi-ordered sets $S$ and $T$ is
called \emph{increasing} if $s\leq
t\Rightarrow\varphi(s)\leq\varphi(t)$ for all $s,t\in S$, and
\emph{strictly increasing} if $s<t\Rightarrow\varphi(s)<\varphi(t)$
for all $s,t\in S$. We also say that $\varphi\colon S\to T$ is a
\emph{quasi-embedding} if $\varphi(s)\leq\varphi(t) \Rightarrow s\leq
t$ for all $s,t\in S$.

An {\it antichain} of $S$ is a subset $A \subseteq S$ such that $s
\not\leq t$ and $t\not\leq s$ for all $s\not\sim t$ in $A$.  A {\it
  final segment} of a quasi-ordered set $(S,\leq)$ is a subset $F
\subseteq S$ which is closed upwards: $s\leq t \ \wedge \ s\in F
\Rightarrow t\in F$, for all $s,t\in S$.  We can view the set ${\cal
  F}(S)$ of final segments of $S$ as an ordered set, with the ordering
given by reverse inclusion.  Given a subset $M$ of $S$, the set
$\bigl\{t \in S : \exists s \in M \text{ with } s \leq t \bigr\}$ is a
final segment of $S$, the final segment {\it generated by $M$.}  An
{\it initial segment} of $S$ is a subset of $S$ whose complement is a
final segment. An initial segment $I$ of $S$ is {\it proper}\/ if
$I\neq S$. For $a\in S$ we denote by $S^{\leq a}$ the initial segment
consisting of all $s\in S$ with $s\leq a$.

A quasi-ordered set $S$ is said to be {\it well-founded} if there is
no infinite strictly decreasing sequence $s_1 > s_2 > \cdots $ in $S$,
and {\it well-quasi-ordered} if in addition every antichain of $S$ is
finite.
The following
characterization of well-quasi-orderings is classical (see, for
example, \cite{Kruskal}). An infinite sequence $s_1,s_2,\dots$ in $S$
is called {\it good}\/ if $s_i\leq s_j$ for some indices $i<j$, and
{\it bad} otherwise.

\begin{prop}\label{equivquasiorder}
  The following are equivalent, for a quasi-ordered set $S$:
\begin{enumerate}
\item $S$ is well-quasi-ordered.
\item Every infinite sequence in $S$ is good.
\item Every infinite sequence in $S$ contains an infinite increasing
  subsequence.
\item Any final segment of $S$ is finitely generated.
\item $\bigl(\mathcal{F} ( S ),\supseteq\bigr)$ is well-founded
  \textup{(}i.e., the ascending chain condition holds for final segments of $S$\textup{).} \qed
\end{enumerate}
\end{prop}

Let $(S,\leq_S)$ and $(T,\leq_T)$ be quasi-ordered sets. If there
exists an increasing surjection $S\to T$ and $S$ is
well-quasi-ordered, then $T$ is well-quasi-ordered, and if there
exists a quasi-embedding $S\to T$ and $T$ is well-quasi-ordered, then
so is $S$.  Moreover, the cartesian product $S\times T$ can be turned
into a quasi-ordered set by using the cartesian product of $\leq_S$ and
$\leq_T$: $$(s,t) \leq (s',t') \quad:\Longleftrightarrow\quad s\leq_S
s' \wedge t\leq_T t', \qquad \text{for $s,s'\in S$, $t,t'\in T$.}$$
Using Proposition~\ref{equivquasiorder} we see that the cartesian
product of two well-quasi-ordered sets is again well-quasi-ordered.

Of course, a total ordering $\leq$ is well-quasi-ordered if and only if
it is well-founded; in this case $\leq$ is called a {\it well-ordering}.
Every well-ordered set
is isomorphic to a unique ordinal number, called its {\it order type}.
The order type of $\N=\{0,1,2,\dots\}$ with its usual ordering is $\omega$.

\subsection{A lemma of Higman}
Given a set $X$, we let $X^*$ denote the set of all finite sequences of
elements of $X$ (including the empty sequence).  We may think of the
elements of $X^*$ as {\em non-commutative words}\/ $x_1\cdots x_m$
with letters $x_1,\dots,x_m$ coming from the alphabet $X$.  With the
concatenation of such words as the operation, $X^*$ is the free monoid
generated by $X$.  A quasi-ordering $\leq$ on $X$ yields a
quasi-ordering $\leq_\Higman$ (the {\em Higman quasi-ordering}\/) on
$X^*$ as follows:
$$x_1\cdots x_m \leq_\Higman y_1\cdots y_n \quad :\Longleftrightarrow
\quad
\begin{cases}
  &\text{\parbox{175pt}{there exists a strictly increasing function
      $\varphi\colon\{1,\dots,m\}\to\{1,\dots,n\}$ such that $x_i \leq
      y_{\varphi(i)}$ for all $1\leq i\leq m$.}}\end{cases}
$$
If $\leq$ is an ordering on $X$, then $\leq_\Higman$ is an ordering
on $X^*$.  The following fact was shown by Higman \cite{Higman} (with
an ingenious proof due to Nash-Williams \cite{NW1}):

\begin{lem}\label{Higman}
  If $\leq$ is a well-quasi-ordering on $X$, then $\leq_\Higman$ is a
  well-quasi-ordering on $X^*$. \qed
\end{lem}

It follows that if $\leq$ is a well-quasi-ordering on $X$, then the
quasi-ordering $\leq^{*}$ on $X^*$ defined by
$$
x_1\cdots x_m \leq^{\ast} y_1\cdots y_n \quad :\Longleftrightarrow
\quad
\begin{cases}
  &\text{\parbox{150pt}{there exists an injective function
      $\varphi\colon \{1,\dots,m\}\to\{1,\dots,n\}$ such that $x_i
      \leq y_{\varphi(i)}$ for all $1\leq i\leq m$}}\end{cases}
$$
is also a well-quasi-ordering (since $\leq^{*}$ extends
$\leq_\Higman$).

We also let $X^\diamond$ be the set of {\em commutative words}\/ in
the alphabet $X$, that is, the free commutative monoid generated by
$X$ (with identity element denoted by $1$).  We sometimes also refer
to the elements of $X^\diamond$ as {\em monomials}\/ (in the set of
indeterminates $X$).  We have a natural surjective monoid homomorphism
$\pi\colon X^*\to X^\diamond$ given by simply ``making the
indeterminates commute'' (i.e., interpreting a non-commutative word
from $X^*$ as a commutative word in $X^\diamond$).  Unlike
$\leq_\Higman$, the quasi-ordering $\leq^{*}$ is compatible with $\pi$
in the sense that $v\leq^{*} w\Rightarrow v'\leq^{*} w'$ for all
$v,v',w,w'\in X^*$ with $\pi(v)=\pi(v')$ and $\pi(w)=\pi(w')$. Hence
$\pi(v)\leq^\diamond \pi(w) :\Longleftrightarrow v\leq^* w$ defines a
quasi-ordering $\leq^\diamond$ on $X^\diamond=\pi(X^*)$ making $\pi$ an
increasing map.  The quasi-ordering $\leq^\diamond$ extends the divisibility
relation in the monoid $X^\diamond$: $$v|w\quad :\Longleftrightarrow
\quad \text{$uv=w$ for some $u\in X^\diamond$.}$$
If we take for
$\leq$ the trivial ordering on $X$, then $\leq^\diamond$ corresponds
exactly to divisibility in $X^\diamond$, and this ordering is a
well-quasi-ordering if and only if $X$ is finite.  In general we have,
as an immediate consequence of Higman's lemma (since $\pi$ is a surjection):

\begin{cor}\label{Higman-Cor}
  If $\leq$ is a well-quasi-ordering on the set $X$, then
  $\leq^\diamond$ is a well-quasi-ordering on $X^\diamond$. \qed
\end{cor}

\subsection{A theorem of Nash-Williams}
Given a totally ordered set $S$ and a quasi-ordered set $X$, we denote
by $\Fin(S,X)$ the set of all functions $f\colon I\to X$, where $I$ is
a proper initial segment of $S$, whose range $f(I)$ is {\it finite.}\/
We define a quasi-ordering $\leq_\Higman$ on $\Fin(S,X)$ as follows:
for $f\colon I\to X$ and $g\colon J\to X$ from $\Fin(S,X)$ put
$$
f \leq_\Higman g \quad :\Longleftrightarrow \quad
\begin{cases}
  &\text{\parbox{220pt}{there exists a strictly increasing function
      $\varphi\colon I\to J$ such that $f(i) \leq g(\varphi(i))$ for
      all $i\in I$.}}\end{cases}
$$
We may think of an element of $\Fin(S,X)$ as a sequence of elements
of $X$ indexed by indices in some proper initial segment of $S$.  So
for $S=\N$ with its usual ordering, we can identify
elements of $\Fin(\N,X)$ with words in $X^*$, and then $\leq_\Higman$
for $\Fin(\N,X)$ agrees with $\leq_\Higman$ on $X^*$ as defined above.
We will have occasion to use a far-reaching generalization of
Lemma~\ref{Higman}:

\begin{thm}\label{NW-Theorem}
  If $X$ is well-quasi-ordered and $S$ is well-ordered, then
  $\Fin(S,X)$ is well-quasi-ordered. \qed
\end{thm}

This theorem was proved by Nash-Williams \cite{NW}; special cases were
shown earlier in \cite{ER, Milner, Rado}.

\subsection{Term orderings}
A {\em term ordering}\/ of $X^\diamond$ is a well-ordering $\leq$ of
$X^\diamond$ such that
\begin{enumerate}
\item $1\leq x$ for all $x\in X$, and
\item $v\leq w\Rightarrow xv\leq xw$ for all $v,w\in X^\diamond$ and
  $x\in X$.
\end{enumerate}
Every ordering $\leq$ of $X^\diamond$ satisfying (1) and (2) extends
the ordering $\leq^\diamond$ obtained from the restriction of $\leq$
to $X$. In particular, $\leq$ extends the divisibility ordering on
$X^\diamond$. By the corollary above, a total ordering $\leq$ of
$X^\diamond$ which satisfies (1) and (2) is a term ordering if and
only if its restriction to $X$ is a well-ordering.

\begin{ex}
  Let $\leq$ be a total ordering of $X$. We define the induced {\em
  lexicographic ordering $\leq_{\lex}$}\/ of monomials as
  follows: given $v,w\in X^\diamond$ we can write $v=x_1^{a_1}\cdots
  x_n^{a_n}$ and $w=x_1^{b_1}\cdots x_n^{b_n}$ with $x_1<\cdots<x_n$
  in $X$ and all $a_i,b_i\in\N$; then
  $$v \leq_{\text{lex}} w \quad :\Longleftrightarrow\quad
  \text{$(a_n,\dots,a_1)\leq (b_n,\dots,b_1)$ lexicographically (from the left).}
  $$
  The ordering $\leq_{\lex}$ is total and satisfies (1), (2);
  hence if the ordering $\leq$ of $X$ is a well-ordering, then
  $\leq_{\lex}$ is a term ordering of $X^\diamond$.
\end{ex}

\begin{rem}
  Let $\leq$ be a total ordering of $X$.
  For $w\in X^\diamond$, $w\neq 1$, we let
  $$|w| := \max\, \{ x\in X: x|w \}\quad \text{(with respect to $\leq$).}$$
  We also put $|1|:=-\infty$,
  where we set $-\infty<x$ for all $x\in X$.  One of the perks of
  using the lexicographic ordering as a term ordering on
  $X^\diamond$ is that if $v$
  and $w$ are monomials with $v\leq_{\lex} w$, then $|v| \leq |w|$. Below,
  we often use this observation.
\end{rem}

The previous example shows that for every set $X$ there exists a term
ordering of $X^\diamond$, since every set can be well-ordered by the
Axiom of Choice.  In fact, every set $X$ can be equipped with a
well-ordering, every proper initial segment of which has strictly
smaller cardinality than $X$; in other words, the order type of this
ordering (a certain ordinal number) is a cardinal number. We shall
call such an ordering of $X$ a {\em cardinal well-ordering}\/ of $X$.

\begin{lem}\label{Extension-Lemma}
  Let $X$ be a set equipped with a cardinal well-ordering, and let $I$
  be a proper initial segment of $X$. Then every injective function
  $I\to X$ can be extended to a permutation of $X$.
\end{lem}
\begin{proof}
  Since this is clear if $X$ is finite, suppose that $X$ is infinite.
  Let $\varphi\colon I\to X$ be injective. Since $I$ has cardinality
  $|I|<|X|$ and $X$ is infinite, we have $|X|=\max\,\{|X\setminus
  I|,|I|\}=|X\setminus I|$.  Similarly, since $|\varphi(I)|=|I|<|X|$,
  we also have $|X\setminus\varphi(I)|=|X|$.  Hence there exists a
  bijection $\psi\colon X\setminus I\to X\setminus\varphi(I)$.
  Combining
  $\varphi$ and $\psi$ yields a permutation of $X$ as desired.
\end{proof}

\subsection{A new ordering of monomials}
Let $G$ be a permutation group on a set $X$, that is, a group $G$
together with a faithful action $(\sigma,x)\mapsto\sigma x\colon
G\times X \to X$ of $G$ on $X$. The action of $G$ on $X$ extends in a
natural way to a faithful action of $G$ on $X^\diamond$: $\sigma w=\sigma x_1\cdots\sigma x_n$ for $\sigma\in G$,
$w=x_1\cdots x_n\in X^\diamond$.  Given a term
ordering $\leq$ of $X^\diamond$, we define a new relation on
$X^\diamond$ as follows:

\begin{defn}\label{defpartialord}(The symmetric cancellation ordering corresponding to $G$ and $\leq$).
\[v \preceq w \quad
:\Longleftrightarrow \quad \begin{cases} &\text{\parbox{200pt}{$v \leq
      w$ and there exist $\sigma \in G$ and a monomial $u \in
      X^\diamond$ such that $w = u \sigma v$ and for all $v' \leq v$,
      we have $u \sigma v' \leq w$.}}\end{cases}\]
\end{defn}

\begin{rem} Every term ordering $\leq$ is \textit{linear}: $v \leq w
  \Longleftrightarrow uv \leq uw$ for all monomials $u,v,w$. Hence the
  condition above may be rewritten as: $v\leq w$ and there exists
  $\sigma\in G$ such that $\sigma v|w$ and $\sigma v' \leq \sigma v$
  for all $v'\leq v$. (We say that ``$\sigma$ witnesses $v\preceq
  w$.'')
\end{rem}

\begin{ex} \label{Main Example}
  Let $X=\{x_1,x_2,\dots\}$ be a countably infinite set of
  indeterminates, ordered such that $x_1<x_2<\cdots$, and let
  ${\leq}={\leq_{\lex}}$ be the corresponding lexicographic
  ordering of $X^\diamond$. Also let $G$ be the group of permutations
  of $\{1,2,3,\dots\}$, acting on $X$ via $\sigma x_i=x_{\sigma(i)}$.
  As an example of the relation $\preceq$, consider the following
  chain:
  \[x_1^2 \preceq x_1x_2^2 \preceq x_1^3x_2x_3^2.\] To verify the
  first inequality, notice that $x_1x_2^2 = x_1 \sigma (x_1^2)$, in
  which $\sigma$ is the transposition $(1\, 2)$.  If $v' =
  x_1^{a_1}\cdots x_n^{a_n} \leq x_1^2$ with $a_1,\dots,a_n\in\N$,
  $a_n>0$, then it follows that $n = 1$ and $a_1 \leq 2$.  In
  particular, $x_1 \sigma v' = x_1x_2^{a_1} \leq x_1x_2^2$.  For the
  second relationship, we have that $x_1^3x_2x_3^2 = x_1^3 \tau
  (x_1x_2^2)$, in which $\tau$ is the cycle $(1\,2\,3)$.
  Additionally, if $v' = x_1^{a_1}\cdots x_n^{a_n} \leq x_1x_2^2$ with
  $a_1,\dots,a_n\in\N$, $a_n>0$, then $n\leq 2$, and if $n = 2$, then
  either $a_2 = 1$ or $a_2=2$, $a_1\leq 1$.  In each case we get
  $x_1^3 \tau v' = x_1^3 x_2^{a_1}x_3^{a_2} \leq x_1^3x_2x_3^2$.
\end{ex}

Although Definition \ref{defpartialord} appears technical, we will soon present a
nice interpretation of it that involves leading term cancellation of
polynomials.  First we verify that it is indeed an ordering.

\begin{lem}
  The relation $\preceq$ is an ordering on monomials.
\end{lem}

\begin{proof}
  First notice that $w \preceq w$ since we may take $u = 1$ and
  $\sigma = $ the identity permutation.  Next, suppose that $u \preceq
  v \preceq w$. Then there exist permutations $\sigma$, $\tau$ in $G$
  and monomials $u_1$, $u_2$ in $X^\diamond$ such that $v = u_1 \sigma
  u$, $w = u_2 \tau v$.  In particular, $w = u_2 (\tau u_1)(\tau
  \sigma u)$.  Additionally, if $v' \leq u$, then $u_1 \sigma v' \leq
  v$, so that $u_2 \tau (u_1 \sigma v') \leq w$.  It follows that $u_2
  (\tau u_1)(\tau \sigma v') \leq w$.  This shows transitivity;
  anti-symmetry of $\preceq$ follows from anti-symmetry of $\leq$.
\end{proof}

We offer a useful interpretation of this ordering (which motivates its
name).  We fix a commutative ring $A$ and let $R=A[X]$ be the ring of
polynomials with coefficients from $A$ in the collection of commuting
indeterminates $X$. Its elements may be written uniquely in the form
$$f=\sum_{w\in X^\diamond} a_w w,$$
where $a_w\in A$ for all $w\in
X^\diamond$, and all but finitely many $a_w$ are zero.  We say that a
monomial $w$ \emph{occurs}\/ in $f$ if $a_w\neq 0$.  Given a non-zero
$f\in R$ we define $\lm(f)$, the {\em leading monomial}\/ of $f$ (with
respect to our choice of term ordering $\leq$) to be the largest
monomial $w$ (with respect to $\leq$) which occurs in $f$.  If
$w=\lm(f)$, then $a_w$ is the \emph{leading coefficient}\/ of $f$,
denoted by $\lc(f)$, and $a_w w$ is the {\em leading term}\/ of $f$,
denoted by $\lt(f)$.  By convention, we set $\lm(0)=\lc(0)=\lt(0)=0$.
We let $R[G]$ be the group ring of $G$ over $R$ (with multiplication
given by $f\sigma\cdot g\tau = fg(\sigma\tau)$ for $f,g\in R$,
$\sigma,\tau\in G$), and we view $R$ as a left $R[G]$-module in the
natural way.

\begin{lem}\label{cancellation}
  Let $f\in R$, $f\neq 0$, and $w\in X^\diamond$. Suppose that
  $\sigma\in G$ witnesses $\lm(f)\preceq w$, and let $u\in X^\diamond$
  with $u\sigma\lm(f)=w$. Then $\lm(u\sigma f)=u\sigma\lm(f)$.
\end{lem}
\begin{proof}
  Put $v=\lm(f)$. Every monomial occurring in $u\sigma f$ has the form
  $u\sigma v'$, where $v'$ occurs in $f$. Hence $v'\leq v$, and since
  $\sigma$ witnesses $v\preceq w$, this yields $u\sigma v'\leq w$.
\end{proof}

Suppose that $A$ is a field, let $v \preceq w$ be in $X^\diamond$ and
let $f$, $g$ be two polynomials in $R$ with leading monomials $v$, $w$,
respectively.  Then, from the definition and
the lemma above, there exists a $\sigma\in G$ and a term $cu$
($c\in A\setminus\{0\}$, $u\in X^\diamond$) such that all monomials
occurring in
\[ h = g - cu
\sigma f\] are strictly smaller (with respect to $\leq$) than $w$.
For readers familiar with the theory of Gr\"obner bases, the
polynomial $h$ can be viewed as a kind of symmetric version of the
$S$-polynomial (see, for instance, \cite[Chapter~15]{eisenbud}).

\begin{ex}
  In the situation of Example~\ref{Main Example} above, let $f =
  x_1x_2^2 + x_2 + x_1^2$ and $g = x_1^3x_2x_3^2 +x_3^2 + x_1^4x_3$.
  Set $\sigma = (1\,2\,3)$, and observe that \[ g - x_1^3 \sigma f =
  x_1^4x_3 + x_3^2 - x_1^3x_3 -x_1^3x_2^2\] has a smaller leading
  monomial than $g$.
\end{ex}

We are mostly interested in the case where our term ordering on
$X^\diamond$ is $\leq_\lex$, and $G=\frak S_X$. Under these
assumptions we have:

\begin{lem}\label{image of sigma}
  Let $v,w\in X^\diamond$ with $v\preceq w$.  Then for every
  $\sigma\in\frak S_X$ witnessing $v\preceq w$ we have $\sigma(X^{\leq
    |v|}) \subseteq X^{\leq |w|}$.  Moreover, if the order type of
  $(X,{\leq})$ is $\leq\omega$, then we can choose such $\sigma$ with
  the additional property that $\sigma(x)=x$ for all $x>|w|$.
\end{lem}
\begin{proof}
  To see the first claim, suppose for a contradiction that $\sigma
  x>|w|$ for some $x\in X$, $x\leq |v|$. We have $\sigma v|w$, so if $x|v$, then
  $\sigma x|w$, contradicting $\sigma x>|w|$. In particular $x<|v|$,
  which yields $x <_\lex v$ and thus $\sigma x \leq_\lex \sigma
  v\leq_\lex w$, again contradicting $\sigma x>|w|$.  Now suppose
  that the order type of $X$ is $\leq\omega$, and let $\sigma$ witness
  $v\preceq w$. Then $|v|\leq |w|$, and $\sigma\upharpoonright X^{\leq
    |v|}$ can be extended to a permutation $\sigma'$ of the finite set
  $X^{\leq |w|}$. We further extend $\sigma'$ to a permutation of $X$
  by setting $\sigma'(x)=x$ for all $x>|w|$. One checks easily that
  $\sigma'$ still witnesses $v\preceq w$.
\end{proof}

\subsection{Lovely orderings}
We say that a term ordering $\leq$ of $X^\diamond$ is {\it lovely}\/
for $G$ if the corresponding symmetric cancellation ordering $\preceq$
on $X^\diamond$ is a well-quasi-ordering.
If $\leq$ is lovely for a subgroup of $G$, then $\leq$ is
lovely for $G$.

\begin{ex}\label{trivial-example}
  The symmetric cancellation ordering corresponding to $G=\{1\}$ and a
  given term ordering $\leq$ of $X^\diamond$ is just
  $$v\preceq w \quad\Longleftrightarrow\quad v\leq w \ \wedge\ v|w.$$
  Hence a term ordering of $X^\diamond$ is lovely for $G=\{1\}$ if and
  only if divisibility in $X^\diamond$ has no infinite antichains;
  that is, exactly if $X$ is finite.
\end{ex}

This terminology is inspired by the following definition from
\cite{Camina-Evans} (which in turn goes back to an idea in
\cite{Ahlbrandt-Ziegler}):

\begin{defn}
  Given an ordering $\leq$ of $X$, consider the following ordering of
  $X$:
\[x \sqsubseteq y \quad
:\Longleftrightarrow \quad \begin{cases} &\text{\parbox{200pt}{$x \leq
      y$ and there exists $\sigma \in G$ such that $\sigma x=y$ and
      for all $x' \leq x$, we have $\sigma x' \leq y$.}}\end{cases}\]
A well-ordering $\leq$ of $X$ is called \emph{nice} (for $G$) if
$\sqsubseteq$ is a well-quasi-ordering.
\end{defn}

In \cite{Ahlbrandt-Ziegler} one finds various examples of nice
orderings, and in \cite{Camina-Evans} it is shown that if $X$ admits a
nice ordering with respect to $G$, then for every field $F$, the free
$F$-module $FX$ with basis $X$ is Noetherian as a module over $F[G]$.
It is clear that the restriction to $X$ of a lovely ordering of
$X^\diamond$ is nice. However, there do exist permutation groups
$(G,X)$ for which $X$ admits a nice ordering, but $X^\diamond$ does
not admit a lovely ordering; see Example~\ref{bad-ex} and
Proposition~\ref{counterexample} below.

\begin{ex}
  Suppose that $X$ is countable. Then every well-ordering of $X$ of
  order type $\omega$ is nice for $\frak S_X$.  To see this, we may
  assume that $X=\N$ with its usual ordering.  It is then easy to see
  that if $x\leq y$ in $\N$, then $x\sqsubseteq y$, witnessed by any
  extension $\sigma$ of the strictly increasing map $n\mapsto
  n+y-x\colon \N^{\leq x}\to\N$ to a permutation of $\N$.
\end{ex}

The following crucial fact (generalizing the last example) is needed
for our proof of Theorem~\ref{onevarfinitegenthm}:

\begin{thm}\label{wellquasithm}
  The lexicographic ordering of $X^\diamond$ corresponding to
  a cardinal well-ordering of a set $X$ is lovely for the full
  symmetric group $\frak S_X$ of $X$.
\end{thm}

For the proof, let as above $\Fin(X,\N)$ be the set of all sequences
in $\N$ indexed by elements in some proper initial segment of $X$
which have finite range, quasi-ordered by $\leq_{\Higman}$.  For a monomial
$w\neq 1$ we define $w^*\colon X^{\leq |w|}\to\N$ by
$$w^*(x) := \max\,\{a\in\N : x^a|w\}.$$
Then clearly
$w^*\in\Fin(X,\N)$; in fact, $w^*(x)=0$ for all but finitely many
$x\in X^{\leq |w|}$. We also let $1^*:=$ the empty sequence
$\emptyset\to\N$ (the unique smallest element of $\Fin(X,\N)$).  We
now quasi-order $X^\diamond\times\Fin(X,\N)$ by the cartesian product of the
ordering $\leq_{\lex}$ on $X^\diamond$ and the quasi-ordering
$\leq_\Higman$ on $\Fin(X,\N)$.  By Corollary~\ref{Higman-Cor},
Theorem~\ref{NW-Theorem}, and the remark following
Proposition~\ref{equivquasiorder}, $X^\diamond\times\Fin(X,\N)$ is
well-quasi-ordered. Therefore, in order to finish the proof of
Theorem~\ref{wellquasithm}, it suffices to show:

\begin{lem}
  The map
  $$w\mapsto (w,w^*)\colon X^\diamond \to X^\diamond\times\Fin(X,\N)$$
  is a quasi-embedding with respect to the symmetric cancellation
  ordering on $X^\diamond$ and the quasi-ordering on
  $X^\diamond\times\Fin(X,\N)$.
\end{lem}
\begin{proof}
  Suppose that $v$, $w$ are monomials with $v\leq_{\lex} w$ and $v^*
  \leq_\Higman w^*$; we need to show that $v\preceq w$. For this we
  may assume that $v,w\neq 1$. So there exists a strictly increasing
  function $\varphi\colon X^{\leq |v|} \to X^{\leq |w|}$ such that
\begin{equation}\label{Div}
v^*(x) \leq w^*(\varphi(x))\qquad\text{ for all $x\in X$ with $x\leq |v|$.}
\end{equation}
By Lemma~\ref{Extension-Lemma} there exists $\sigma\in\frak S_X$ such
that $\sigma \upharpoonright X^{\leq |v|}=\varphi \upharpoonright
X^{\leq |v|}$.  Then clearly $\sigma v|w$ by \eqref{Div}. Now let
$v'\leq_{\lex} v$; we claim that $\sigma v'\leq_\lex \sigma v$.
Again we may assume $v'\neq 1$.  Then $|v'|\leq |v|$; hence we may
write $$v'=x_1^{a_1}\cdots x_n^{a_n}, \quad v=x_1^{b_1}\cdots
x_n^{b_n}$$
with $x_1<\cdots< x_n\leq |v|$ in $X$ and $a_i,b_j\in\N$.
Put $y_1:=\varphi(x_1),\dots,y_n:=\varphi(x_n)$.  Then
$y_1<\cdots<y_n$ and
$$\sigma v'=y_1^{a_1}\cdots y_n^{a_n},\quad \sigma v=y_1^{b_1}\cdots
y_n^{b_n},$$
and therefore $\sigma v' \leq_\lex \sigma v$ as
required.
\end{proof}

\subsection{The case of countable $X$}

In Section~4 we will apply Theorem~\ref{wellquasithm} in the case
where $X$ is countable. Then the order type of $X$ is at most
$\omega$, and in the proof of the theorem given above we only need to
appeal to a special instance (Higman's Lemma) of
Theorem~\ref{NW-Theorem}. We finish this section by giving a
self-contained proof of this important special case of
Theorem~\ref{wellquasithm}, avoiding Theorem~\ref{NW-Theorem}.  Let
$\frak S_{(X)}$ denote the subgroup of $\frak S_X$ consisting of all
$\sigma\in\frak S_X$ with the property that $\sigma(x)=x$ for all but
finitely many letters $x\in X$.

\begin{thm}\label{wellquasithm2}
  The lexicographic ordering of $X^\diamond$ corresponding to
  a cardinal well-ordering of a countable set $X$ is lovely for $\frak
  S_{(X)}$.
\end{thm}

Let $X$ be countable and let $\leq$ be a cardinal well-ordering of $X$.
Enumerate the elements of
$X$ as $x_1<x_2<\cdots$.
We assume that $X$ is infinite;
this is not a restriction, since
by Lemma~\ref{image of sigma} we have:

\begin{lem}
  If the lexicographic ordering of $X^\diamond$ is lovely for
  $\frak S_{(X)}$, then for any $n$ and $X_n:=\{x_1,\dots,x_n\}$, the
  lexicographic ordering of $(X_n)^\diamond$ is lovely for
  $\frak S_{X_n}$. \qed
\end{lem}

We begin with some preliminary lemmas. Here, $\preceq$ is the
symmetric cancellation ordering corresponding to $\frak S_{(X)}$ and
$\leq_\lex$.  We identify $\frak S_{(X)}$ and $\frak S_\infty:=\frak
S_{(\N)}$ in the natural way, and for every $n$ we regard $\frak S_n$,
the group of permutations of $\{1,2,\dots,n\}$, as a subgroup of
$\frak S_{\infty}$; then $\frak S_n\leq \frak S_{n+1}$ for each $n$,
and $\frak S_{\infty}=\bigcup_n \frak S_n$.

\begin{lem}\label{oneshiftuplem}
  Suppose that $x_1^{a_1}\cdots x_n^{a_n} \preceq x_1^{b_1}\cdots
  x_n^{b_n}$, where $a_i,b_j\in\N$, $b_n>0$. Then for any $c \in \N$ we
  have $x_1^{a_1}\cdots x_n^{a_n} \preceq x_1^c x_2^{b_1}\cdots
  x_{n+1}^{b_n}$.
\end{lem}

\begin{proof}
  Let $v:=x_1^{a_1}\cdots x_n^{a_n}$, $w:= x_1^{b_1}\cdots x_n^{b_n}$.
  We may assume $v\neq 1$.  Clearly $v \leq_\lex w$ and $b_n>0$ yield
  $x_1^{a_1}\cdots x_n^{a_n} \leq_\lex x_1^c x_2^{b_1}\cdots
  x_{n+1}^{b_n}$.  Now let $\sigma\in\frak S_\infty$ witness $v\preceq
  w$.  Let $\tau$ be the cyclic permutation $\tau = (1\,2\,3 \cdots
  (n+1))$ and set $\hat\sigma:=\tau\sigma$.  Then $\sigma v|w$ yields
  $\hat\sigma v|\tau w$, hence $\hat\sigma v|x_1^c\tau w = x_1^c
  x_2^{b_1}\cdots x_{n+1}^{b_n}$.  Next, suppose that $v'\leq_\lex
  v$; then $\sigma v'\leq_\lex \sigma v$.  By Lemma~\ref{image of
    sigma} and the nature of $\tau$, the map $\tau\upharpoonright \sigma(\{1,\dots,|v|\})$ is strictly
  increasing, which gives $\hat\sigma v'=\tau\sigma
  v'\leq_{\lex}\tau\sigma v=\hat\sigma v$.  Hence $\hat\sigma$
  witnesses $x_1^{a_1}\cdots x_n^{a_n} \preceq x_1^c x_2^{b_1}\cdots
  x_{n+1}^{b_n}$.
\end{proof}

\begin{lem}\label{twoshiftuplem}
  If $x_1^{a_1}\cdots x_n^{a_n} \preceq x_1^{b_1}\cdots x_n^{b_n}$,
  where $a_i,b_j\in\N$, $b_n>0$, and $a,b \in \N$ are such that $a
  \leq b$, then $x_1^ax_2^{a_1}\cdots x_{n+1}^{a_n} \preceq
  x_1^bx_2^{b_1}\cdots x_{n+1}^{b_{n}}$.
\end{lem}

\begin{proof}
  As before let $v:=x_1^{a_1}\cdots x_n^{a_n}$, $w:= x_1^{b_1}\cdots
  x_n^{b_n}$. Once again, we may assume $v\neq 1$, and it is clear
  that $x_1^ax_2^{a_1}\cdots x_{n+1}^{a_n} \leq_\lex
  x_1^bx_2^{b_1}\cdots x_{n+1}^{b_{n}}$.  Let $\sigma\in\frak
  S_\infty$ witness $v\preceq w$.  By Lemma~\ref{image of sigma} we
  may assume that $\sigma(x_i)=x_i$ for all $i>n$.  Let $\tau$ be the
  cyclic permutation $\tau = (1\,2 \cdots (n+1))$.  Setting
  $\hat{\sigma} = \tau \sigma \tau^{-1}$, we have $\hat\sigma x_1 =
  x_1$; hence
\begin{equation}\label{sigma hat}
  \hat{\sigma} (x_1^a x_{2}^{a_1} \cdots x_{n+1}^{a_n} ) =
  \hat\sigma(x_1^a) \hat\sigma(x_{2}^{a_1} \cdots x_{n+1}^{a_n} )
   =   x_{1}^a \tau\sigma v.
\end{equation}
Since $\sigma v|w$, this last expression divides $x_1^b\tau w =
x_{1}^{b}x_{2}^{b_1} \cdots x_{n+1}^{b_n}$.  Suppose that $v' =
x_1^{c_1} \cdots x_{n+1}^{c_{n+1}} \leq_\lex x_1^{a}x_{2}^{a_1}
\cdots x_{n+1}^{a_n}$, where $c_i\in\N$.  Then, since we are using
a lexicographic order, we have \[x_2^{c_2} \cdots
x_{n+1}^{c_{n+1}} \leq_\lex x_{2}^{a_1} \cdots x_{n+1}^{a_n}\] and
therefore
\[ \tau^{-1}(x_2^{c_2} \cdots x_{n+1}^{c_{n+1}})=x_1^{c_2} \cdots x_{n}^{c_{n+1}}
\leq_\lex \tau^{-1}(x_{2}^{a_1} \cdots x_{n+1}^{a_n})=v.\] By
assumption, this implies that $\sigma \tau^{-1} (x_2^{c_2} \cdots
x_{n+1}^{c_{n+1}}) \leq_\lex \sigma v$ and thus by \eqref{sigma hat},
$$\hat{\sigma} (x_2^{c_2} \cdots x_{n+1}^{c_{n+1}}) \leq_\lex
\tau\sigma v=\hat\sigma(x_{2}^{a_1} \cdots x_{n+1}^{a_n} ).$$
If this
inequality is strict, then since $1\notin\hat\sigma\big(\{2,\dots,n+1\}\big)$, clearly
$$\hat\sigma v'=x_1^{c_1} \hat{\sigma} (x_2^{c_2} \cdots
x_{n+1}^{c_{n+1}})<_\lex x_1^a\tau\sigma v= \hat{\sigma} (x_1^a
x_{2}^{a_1} \cdots x_{n+1}^{a_n} ).$$
Otherwise $x_2^{c_2} \cdots
x_{n+1}^{c_{n+1}}=x_{2}^{a_1} \cdots x_{n+1}^{a_n}$; hence $c_1 \leq
a$, in which case we still have $\hat{\sigma} v' \leq_\lex
\hat{\sigma}(x_1^a x_{2}^{a_1} \cdots x_{n+1}^{a_n} )$.  Therefore
$\hat\sigma$ witnesses $x_1^ax_2^{a_1}\cdots x_{n+1}^{a_n} \preceq
x_1^bx_2^{b_1}\cdots x_{n+1}^{b_{n}}$.  This completes the proof.
\end{proof}

We now have enough to show Theorem~\ref{wellquasithm2}.  The proof
uses the basic idea from Nash-Williams' proof \cite{NW} of Higman's
lemma. Assume for the sake of contradiction that there exists a bad
sequence $$w^{(1)},w^{(2)},\dots,w^{(n)},\dots \qquad\text{in
  $X^\diamond$.}$$
For $w\in X^\diamond\setminus\{1\}$ let $j(w)$ be
the index $j\geq 1$ with $|w|=x_j$, and put $j(1):=0$. We may assume
that the bad sequence is chosen in such a way that for every $n$,
$j(w^{(n)})$ is {\it minimal} among the $j(w)$, where $w$ ranges over
all elements of $X^\diamond$ with the property that
$w^{(1)},w^{(2)},\dots,w^{(n-1)},w$ can be continued to a bad sequence
in $X^\diamond$.  Because $1\leq_{\lex} w$ for all $w\in X^\diamond$, we have $j(w^{(n)})>0$ for all $n$. For
every $n>0$, write $w^{(n)}=x_1^{a^{(n)}}v^{(n)}$ with $a^{(n)}\in\N$
and $v^{(n)}\in X^\diamond$ not divisible by $x_1$.  Since $\N$ is
well-ordered, there is an infinite sequence $1\leq i_1<i_2<\cdots$ of
indices such that $a^{(i_1)} \leq a^{(i_2)} \leq \cdots$.  Consider
the monoid homomorphism $\alpha\colon X^\diamond\to X^\diamond$ given
by $\alpha(x_{i+1})=x_i$ for all $i>1$. Then $j(\alpha(w))=j(w)-1$ if
$w\neq 1$. Hence by minimality of $w^{(1)},w^{(2)},\dots$, the
sequence
$$w^{(1)},w^{(2)},\dots,w^{(i_1-1)},\alpha(v^{(i_1)}),\alpha(v^{(i_2)}),\dots,\alpha(v^{(i_n)}),\dots
$$
is good; that is, there exist $j < i_1$ and $k$ with $w^{(j)}
\preceq \alpha(v^{(i_k)})$, or there exist $k<l$ with
$\alpha(v^{(i_k)}) \preceq \alpha(v^{(i_l)})$.  In the first case we
have $w^{(j)}\preceq w^{(i_k)}$ by Lemma \ref{oneshiftuplem}; and in
the second case, $w^{(i_k)}\preceq w^{(i_l)}$ by Lemma
\ref{twoshiftuplem}. This contradicts the badness of our sequence
$w^{(1)},w^{(2)},\dots$, finishing the proof.

\begin{question}
  Careful inspection of the proof of Theorem~\ref{wellquasithm} (in
  particular Lemma~\ref{Extension-Lemma}) shows that in the statement
  of the theorem, we can replace $\frak S_X$ by its subgroup
  consisting of all $\sigma$ with the property that the set of $x\in
  X$ with $\sigma(x)\neq x$ has cardinality $<|X|$.  In
  Theorem~\ref{wellquasithm}, can one always replace $\frak S_X$ by
  $\frak S_{(X)}$?
\end{question}

\section{Proof of the Finiteness Theorem}\label{mainthmproof}

We now come to the proof our main result. Throughout this section we
let $A$ be a commutative Noetherian ring, $X$ an arbitrary set,
$R=A[X]$, and we let $G$ be a permutation group on $X$.  An
$R[G]$-submodule of $R$ will be called a \emph{$G$-invariant ideal} of
$R$, or simply an \emph{invariant ideal}, if $G$ is understood.  We
will show:

\begin{thm}\label{finitenessthm}
  If $X^\diamond$ admits a lovely term ordering for $G$, then $R$ is
  Noetherian as an $R[G]$-module.
\end{thm}

For $G=\{1\}$ and $X$ finite, this theorem reduces to Hilbert's basis
theorem, by Example~\ref{trivial-example}.  We also obtain Theorem
\ref{onevarfinitegenthm}:

\begin{cor}
  The $R[\frak S_X]$-module $R$ is Noetherian.
\end{cor}
\begin{proof}
  Choose a cardinal well-ordering of $X$. Then the corresponding
  lexicographic ordering of $X^\diamond$ is lovely for $\frak
  S_X$, by Theorem~\ref{wellquasithm}. Apply
  Theorem~\ref{finitenessthm}.
\end{proof}

\begin{rem}
  It is possible to replace the use of Theorem~\ref{wellquasithm} in the proof
  of the corollary above by the more elementary Theorem~\ref{wellquasithm2}.
  This is because if the $R[\frak S_X]$-module $R$ were not Noetherian, then
  one could find a {\it countably generated}\/
  $R[\frak S_X]$-submodule of $R$ which
  is not finitely generated, and hence a countable subset $X'$ of $X$ such that
  $R'=A[X']$ is not a Noetherian $R'[\frak S_{X'}]$-module.
\end{rem}

The following example shows how the conclusion of
Theorem~\ref{finitenessthm} may fail:

\begin{ex}\label{bad-ex}
  Suppose that $G$ has a cyclic subgroup $H$ which acts freely and
  transitively on $X$. Then $X$ has
  a nice ordering (see \cite{Ahlbrandt-Ziegler}), but
  $R=\Q[X^\diamond]$ is not Noetherian. To see this let $\sigma$ be a
  generator for $H$, and let $x\in X$ be arbitrary. Then the
  $R[G]$-submodule of $R=\Q[X^\diamond]$ generated by the elements
  $\sigma^nx\sigma^{-n}x$ ($n\in\N$) is not finitely generated. So
  by Theorem~\ref{finitenessthm}, $X^\diamond$ does not admit a lovely
  term ordering for $G$.
\end{ex}

For the proof of Theorem~\ref{finitenessthm} we develop a bit of
Gr\"obner basis theory for the $R[G]$-module $R$.  For the time being,
we fix an arbitrary term ordering $\leq$ (not necessarily lovely for
$G$) of $X^\diamond$.

\subsection{Reduction of polynomials}
Let $f\in R$, $f\neq 0$, and let $B$ be a set of non-zero polynomials
in $R$. We say that $f$ is \emph{reducible by $B$}\/ if there exist
pairwise distinct $g_1,\dots,g_m\in B$, $m\geq 1$, such that for each
$i$ we have $\lm(g_i)\preceq \lm(f)$, witnessed by some $\sigma_i\in
G$, and
$$\lt(f) = a_1w_1\sigma_1\lt(g_1) + \cdots + a_mw_m\sigma_m\lt(g_m)$$
for non-zero $a_i\in A$ and monomials $w_i\in X^\diamond$ such that
$w_i\sigma_i\lm(g_i)=\lm(f)$.  In this case we write
$f\underset{B}\longrightarrow h$, where
$$h=f - \big(a_1w_1\sigma_1g_1 + \cdots + a_mw_m\sigma_mg_m\big),$$
and we say that $f$ \emph{reduces to $h$} by $B$.  We say that $f$ is
\emph{reduced} with respect to $B$ if $f$ is not reducible by $B$. By
convention, the zero polynomial is reduced with respect to $B$.
Trivially, every element of $B$ reduces to $0$.

\begin{ex}
  Suppose that $A$ is a field. Then $f$ is reducible by $B$ if and
  only if there exists some $g\in B$ such that $\lm(g)\preceq\lm(f)$.
\end{ex}

\begin{ex}\label{reduction-example}
  Suppose that $f$ is reducible by $B$ as defined (for finite $X$) in,
  say, \cite[Chapter 4]{AL}; that is, there exist $g_1,\dots,g_m\in B$
  and $a_1,\dots,a_m\in A$ ($m\geq 1$) such that $\lm(g_i)|\lm(f)$ for
  all $i$ and
  $$\lc(f) = a_1\lc(g_1) + \cdots + a_m\lc(g_m).$$
  Then $f$ is
  reducible by $B$ in the sense defined above (taking $\sigma_i=1$
  for all $i$).
\end{ex}

\begin{rem}\label{reduction-remark}
  Suppose that $G=\frak S_{X}$, the term ordering $\leq$ of
  $X^\diamond$ is $\leq_{\lex}$, and the order type of $(X,\leq)$ is
  $\leq\omega$. Then in the definition of reducibility by $B$ above,
  we may require that the $\sigma_i$ satisfy $\sigma_i(x)=x$ for all
  $1\leq i\leq m$ and $x>|\lm(f)|$ (by Lemma~\ref{image of sigma}).
\end{rem}

The smallest quasi-ordering on $R$ extending the relation
$\underset{B}\longrightarrow$ is denoted by
$\underset{B}{\overset{*}\longrightarrow}$.  If $f,h\neq 0$ and
$f\underset{B}\longrightarrow h$, then $\lm(h)<\lm(f)$, by
Lemma~\ref{cancellation}.  In particular, every chain
$$h_0\underset{B}\longrightarrow h_1\underset{B}\longrightarrow h_2
\underset{B}\longrightarrow \cdots$$
with all $h_i\in R\setminus\{0\}$
is finite (since the term ordering $\leq$ is well-founded).  Hence
there exists $r\in R$ such that
$f\underset{B}{\overset{*}\longrightarrow} r$ and $r$ is reduced with
respect to $B$; we call such an $r$ a \emph{normal form} of $f$ with
respect to $B$.

\begin{lem}\label{reduction}
  Suppose that $f\underset{B}{\overset{*}\longrightarrow} r$. Then
  there exist $g_1,\dots,g_n\in B$, $\sigma_1,\dots,\sigma_n\in G$ and
  $h_1,\dots,h_n\in R$ such that
  $$f=r+\sum_{i=1}^n h_i\sigma_i g_i\quad \text{and}\quad \lm(f)\geq
  \max_{1\leq i\leq n}\lm(h_i\sigma_ig_i).$$
  \textup{(}In particular,
  $f-r\in \<B\>_{R[G]}$.\textup{)}
\end{lem}
\begin{proof}
  This is clear if $f=r$. Otherwise we have $f
  \underset{B}\longrightarrow h
  \underset{B}{\overset{*}\longrightarrow} r$ for some $h\in R$.
  Inductively we may assume that there exist $g_1,\dots,g_n\in B$,
  $\sigma_1,\dots,\sigma_n\in G$
  and $h_1,\dots,h_n\in R$ such that
  $$h=r+\sum_{i=1}^n h_i\sigma_i g_i\quad \text{and}\quad \lm(h)\geq
  \max_{1\leq i\leq n}\lm(h_i\sigma_ig_i).$$
  There are also
  $g_{n+1},\dots,g_{n+m}\in B$, $\sigma_{n+1},\dots,\sigma_{n+m}\in
  G$, $a_{n+1},\dots,a_{n+m}\in A$ and $w_{n+1},\dots,w_{n+m}\in
  X^\diamond$ such that $\lm(w_{n+i}\sigma_{n+i}g_{n+i})=\lm(f)$ for
  all $i$ and
\[
\lt(f) = \sum_{i=1}^m a_{n+i}w_{n+i}\sigma_{n+i}\lt(g_{n+i}), \qquad f
= h + \sum_{i=1}^m a_{n+i} w_{n+i}\sigma_{n+i} g_{n+i}.\] Hence
putting $h_{n+i}:=a_{n+i}w_{n+i}$ for $i=1,\dots,m$ we have
$f=r+\sum_{j=1}^{n+m} h_j\sigma_j g_j$ and $\lm(f)>\lm(h)\geq
\lm(h_j\sigma_jg_j)$ if $1\leq j\leq n$, $\lm(f)=\lm(h_j\sigma_jg_j)$
if $n<j\leq n+m$.
\end{proof}

\begin{rem}\label{reduction-remark, 2}
  Suppose that $G=\frak S_{X}$, ${\leq}={\leq_\lex}$, and $X$ has
  order type $\leq\omega$.  Then in the previous lemma we can choose
  the $\sigma_i$ such that in addition $\sigma_i(x)=x$ for all $i$ and
  all $x>|\lm(f)|$ (by Remark~\ref{reduction-remark}).
\end{rem}

\subsection{Gr\"obner bases}
Let $B$ be a subset of $R$. We let $$\lt(B):=\big\<\lc(g)w: 0\neq g\in
B,\ \lm(g)\preceq w\big\>_{A}$$
be the $A$-submodule of $R$ generated
by all elements of the form $\lc(g) w$, where $g\in B$ is non-zero and
$w$ is a monomial with $\lm(g)\preceq w$.  Clearly for non-zero $f\in
R$ we have: $\lt(f)\in\lt(B)$ if and only if $f$ is reducible by $B$.
In particular, $\lt(B)$ contains $\big\{\lt(g):g\in B\big\}$, and for
an ideal $I$ of $R$ which is $G$-invariant, we simply have (using
Lemma~\ref{cancellation})
$$\lt(I) = \big \< \lt(f) : f\in I \big\>_A.$$
\begin{defn}
We say that a subset $B$ of an invariant
ideal $I$ of $R$ is a \emph{Gr\"obner basis}\/ for $I$ (with respect
to our choice of term ordering $\leq$) if $\lt(I)=\lt(B)$.
\end{defn}

Additionally, in the case when $A$ is a field,
a Gr\"obner basis is called \textit{minimal} if
no leading monomial of an element in $B$ is $\preceq$ smaller
than any other leading monomial of an element in $B$.


\begin{lem}\label{char GB}
  Let $I$ be an invariant ideal of $R$ and $B$ be a set of non-zero
  elements of $I$. The following are equivalent:
\begin{enumerate}
\item $B$ is a Gr\"obner basis for $I$.
\item Every non-zero $f\in I$ is reducible by $B$.
\item Every $f\in I$ has normal form $0$. \textup{(}In particular,
  $I=\<B\>_{R[G]}$.\textup{)}
\item Every $f\in I$ has unique normal form $0$.
\end{enumerate}
\end{lem}
\begin{proof}
  The implications
  (1)~$\Rightarrow$~(2)~$\Rightarrow$~(3)~$\Rightarrow$~(4) are either
  obvious or follow from the remarks preceding the lemma.  Suppose
  that (4) holds. Every $f\in I\setminus\{0\}$ with
  $\lt(f)\notin\lt(B)$ is reduced with respect to $B$, hence has two
  distinct normal forms ($0$ and $f$), a contradiction. Thus
  $\lt(I)=\lt(B)$.
\end{proof}

Suppose that $B$ is a Gr\"obner basis for an ideal $I$ of the
polynomial ring $R=A[X^\diamond]$, in the usual sense of the word (as
defined, for finite $X$, in \cite[Chapter~4]{AL}); if $I$ is
invariant, then $B$ is a Gr\"obner basis for $I$ as defined above (by
Example~\ref{reduction-example}).  Moreover, for $G=\{1\}$, the
previous lemma reduces to a familiar characterization of Gr\"obner
bases in the usual case of polynomial rings.  It is probably possible
to also introduce a notion of an $S$-polynomial and to prove a
Buchberger-style criterion for Gr\"obner bases in our setting, leading
to a completion procedure for the construction of Gr\"obner bases.  At
this point, we will not pursue these issues further, and rather show:

\begin{prop}\label{finite GB}
  Suppose that the term ordering $\leq$ of $X^\diamond$ is lovely for
  $G$.  Then every invariant ideal of $R$ has a finite Gr\"obner
  basis.
\end{prop}

For a subset $B$ of $R$ let $\lm(B)$ denote the final segment of
$X^\diamond$ with respect to $\preceq$ generated by the $\lm(g)$,
$g\in B$. If $A$ is a field, then a subset $B$ of an invariant ideal
$I$ of $R$ is a Gr\"obner basis for $I$ if and only if
$\lm(B)=\lm(I)$. Hence in this case, the proposition follows
immediately from the equivalence of (1) and (4) in
Proposition~\ref{equivquasiorder}.  For the general
case we use the following observation:

\begin{lem}
  Let $S$ be a well-quasi-ordered set and $T$ be a well-founded
  ordered set, and let $\varphi\colon S\to T$ be decreasing: $s\leq
  t\Rightarrow \varphi(s)\geq\varphi(t)$, for all $s,t\in S$. Then the
  quasi-ordering $\leq_\varphi$ on $S$ defined by
  $$s \leq_\varphi t \quad:\Longleftrightarrow \quad s\leq t\ \wedge\
  \varphi(s)=\varphi(t)$$
  is a well-quasi-ordering. \qed
\end{lem}

\begin{proof}[Proof of Proposition~\ref{finite GB}]
Suppose now that our term ordering of $X^\diamond$ is lovely for $G$,
and let $I$ be an invariant ideal of $R$.  For $w\in X^\diamond$
consider
$$\lc(I,w) := \big\{ \lc(f) : \text{$f\in I$, and $f=0$ or $\lm(f) =
  w$} \big\},$$
an ideal of $A$.  Note that if $v\preceq w$, then
$\lc(I,v)\subseteq\lc(I,w)$.  We apply the lemma to $S=X^\diamond$,
quasi-ordered by $\preceq$, $T=$ the collection of all ideals of $A$,
ordered by reverse inclusion, and $\varphi$ given by $w\mapsto
\lc(I,w)$. Thus by (4) in
Proposition~\ref{equivquasiorder}, applied to the
final segment $X^\diamond$ of the well-quasi-ordering $\leq_\varphi$,
we obtain finitely many $w_1,\dots,w_m\in X^\diamond$
with the following property: for every $w\in X^\diamond$ there exists
some $i\in\{1,\dots,m\}$ such that $w_i\preceq w$ and
$\lc(I,w_i)=\lc(I,w)$. Using Noetherianity of $A$, for every $i$ we
now choose finitely many non-zero elements $g_{i1},\dots,g_{in_i}$ of
$I$ ($n_i\in\N$), each with leading monomial $w_i$, whose leading
coefficients generate the ideal $\lc(I,w_i)$ of $A$.  We claim that
$$B:=\{ g_{ij} : 1\leq i\leq m,\ 1\leq j\leq n_i\}$$
is a Gr\"obner
basis for $I$. To see this, let $0\neq f\in I$, and put $w:=\lm(f)$.
Then there is some $i$ with $w_i\preceq w$ and $\lc(I,w_i)=\lc(I,w)$.
This shows that $f$ is reducible by $\{g_{i1},\dots,g_{i,n_i}\}$, and
hence by $B$. By Lemma~\ref{char GB}, $B$ is a Gr\"obner basis for
$I$.
\end{proof}

From Proposition~\ref{finite GB} and
the implication (1)~$\Rightarrow$~(3) in Lemma~\ref{char GB}
we obtain Theorem~\ref{finitenessthm}.

\subsection{A  partial converse of Theorem~\ref{finitenessthm}}
Consider now the quasi-ordering $|_G$ of $X^\diamond$ defined by
$$v |_G w \quad :\Longleftrightarrow \quad
\exists\sigma\in G: \sigma v|w,$$
which extends every symmetric cancellation ordering corresponding
to a  term ordering of $X^\diamond$.
If $M$ is a set of monomials from $X^\diamond$ and  $F$ the
final segment of $(X^\diamond,{|_G})$
generated by $M$, then the
invariant ideal $\<M\>_{R[G]}$ of $R$ is finitely generated as an
$R[G]$-module if and only if $F$ is generated by a finite subset of $M$.
Hence by the implication (4)~$\Rightarrow$~(1) in
Proposition~\ref{equivquasiorder} we get:

\begin{lem}\label{converse of finitenessthm}
If $R$ is Noetherian as an $R[G]$-module, then $|_G$ is a well-quasi-ordering. \qed
\end{lem}

This will be used in Section~\ref{chemmotivation} below.

\subsection{Connection to a concept due to Michler}

Let $\leq$ be a term ordering of $X^\diamond$.  For each $\sigma\in G$
we define a term ordering $\leq_\sigma$ on $X^\diamond$ by
$$v \leq_\sigma w \quad \Longleftrightarrow \quad \sigma v\leq \sigma
w.$$
We denote the leading monomial of $f\in R$ with respect to
$\leq_\sigma$ by $\lm_\sigma(f)$. Clearly we have
\begin{equation}\label{sigma lm(f)}
\sigma\lm(f) = \lm_{\sigma^{-1}}(\sigma f)\qquad\text{for all $\sigma\in G$ and $f\in R$.}
\end{equation}
Let $I$ be an invariant ideal of $R$.  Generalizing terminology
introduced in \cite{Michler}, let us call a set $B$ of non-zero
elements of $I$ a \emph{universal $G$-Gr\"obner basis} for $I$ (with
respect to $\leq$) if $B$ contains, for every $\sigma\in G$, a
Gr\"obner basis (in the usual sense of the word) for the ideal $I$
with respect to the term ordering $\leq_\sigma$.  If the set $X$ of
indeterminates is finite, then every invariant ideal of $R$ has a
finite universal $G$-Gr\"obner basis.  By the remark following
Lemma~\ref{char GB}, every universal $G$-Gr\"obner basis for an
invariant ideal $I$ of $R$ is a Gr\"obner basis for $I$.  We finish
this section by observing:

\begin{lem}
  Suppose that $A$ is a field. If $B$ is a Gr\"obner basis for the
  invariant ideal $I$ of $R$, then $$GB=\{\sigma g:\sigma\in G,\ g\in
  B\}$$
  is a universal $G$-Gr\"obner basis for $I$.
\end{lem}
\begin{proof}
  Let $\sigma\in G$ and $f\in I$, $f\neq 0$. Then $\sigma f\in I$;
  hence there exists $\tau\in G$ and $g\in B$ such that $w\leq\lm(g)
  \Rightarrow w\leq_\tau\lm(g)$ for all $w\in X^\diamond$, and
  $\tau\lm(g)|\lm(\sigma f)$.  The first condition implies in
  particular that $\tau\lm(g)=\lm(\tau g)$; hence
  $\sigma^{-1}\tau\lm(g)=\lm_\sigma(\sigma^{-1}\tau g)$ and
  $\sigma^{-1}\lm(\sigma f)=\lm_\sigma(f)$ by \eqref{sigma lm(f)}.
  Put $h:=\sigma^{-1}\tau g\in GB$.  Then
  $\lm_\sigma(h)|\lm_\sigma(f)$ by the second condition. This shows
  that $GB$ contains a Gr\"obner basis for $I$ with respect to
  $\leq_\sigma$, as required.
\end{proof}

\begin{ex}
  Suppose that $G=\frak S_n$, the group of permutations of
  $\{1,2,\dots,n\}$, acting on $X=\{x_1,\dots,x_n\}$ via $\sigma
  x_i=x_{\sigma(i)}$.  The invariant ideal $I=\<x_1,\dots,x_n\>_R$ has
  Gr\"obner basis $\{x_1\}$ with respect to the lexicographic
  ordering; a corresponding (minimal) universal $\frak S_n$-Gr\"obner
  basis for $I$ is $\{x_1,\dots,x_n\}$.
\end{ex}

\section{Invariant chains of ideals}

In this section we describe a relationship between certain chains of
increasing ideals in finite-dimensional polynomials rings and
invariant ideals of infinite-dimensional polynomial rings.  We begin
with an abstract setting that is suitable for placing the motivating
problem (described in the next section) in a proper context.
Throughout this section, $m$ and $n$ range over the set of positive
integers. For each $n$, let $R_n$ be a commutative ring, and assume
that $R_n$ is a subring of $R_{n+1}$, for each $n$. Suppose that the
symmetric group on $n$ letters ${\mathfrak S}_n$ gives an action (not
necessarily faithful) on $R_n$  such that $f\mapsto\sigma f\colon
R_n\to R_n$ is a ring homomorphism, for each $\sigma \in {\mathfrak
  S}_n$. Furthermore, suppose that the natural embedding of
${\mathfrak S}_n$ into ${\mathfrak S}_m$ for $n \leq m$ is compatible
with the embedding of rings $R_n\subseteq R_m$; that is, if $\sigma \in {\mathfrak
  S}_n$ and $\hat{{\sigma}}$ is the corresponding element in
${\mathfrak S}_m$, then $ {\hat{{\sigma}} } \upharpoonright {R_n } =
\sigma$.  Note that there exists a unique action of $\frak S_{\infty}$
on the ring $R := \bigcup_{n\geq 1} R_n$ which extends the action of
each $\frak S_n$ on $R_n$. An ideal of $R$ is {\em invariant} if
$\sigma f\in I$ for all $\sigma\in\frak S_\infty$, $f\in I$.

We will need a method for lifting ideals of smaller rings into larger
ones, and one such technique is as follows.

\begin{defn}  For $m \geq n$,
  the \textit{$m$-symmetrization} $L_m(B)$ of a set $B$ of elements of
  $R_n$  
  is the $\frak S_m$-invariant ideal of $R_m$ given by
\[L_m(B) = \< g : g
\in B \>_{R_m[{\mathfrak S}_m]}.\]
\end{defn}

In order for us to apply this definition sensibly, we must make sure
that the $m$-symmetrization of an ideal can be defined in terms of
generators.

\begin{lem}\label{symmetrgeneratorslemma}
  If $B$ is a set of generators for the ideal $I_B = \<B\>_{R_n}$ of
  $R_n$, then $L_m(I_B) = L_m(B)$.
\end{lem}

\begin{proof}
  Suppose that $B$ generates the ideal $I_B \subseteq R_n$. Clearly,
  $L_m(B) \subseteq L_m(I_B)$.  Therefore, it is enough to show the
  inclusion $L_m(I_B) \subseteq L_m(B)$.  Suppose that $h \in
  L_m(I_B)$ so that $h = \sum_{j=1}^{s} f_{j} \cdot \sigma_j h_j$ for
  elements $f_j \in R_m$, $h_j \in I_B$ and $\sigma_j \in {\mathfrak
    S}_m$. Next express each $h_j = \sum_{i=1}^{r_j} p_{ij} g_{ij}$
  for $p_{ij} \in R_n$ and $g_{ij} \in B$.  Substitution into the
  expression above for $h$ gives us
\[ h = \sum_{j=1}^{s}\sum_{i=1}^{r_j} f_{j} \cdot \sigma_j p_{ij}
\cdot \sigma_j g_{ij}.\] This is easily seen to be an element of
$L_m(B)$, completing the proof.
\end{proof}

\begin{ex}\label{naturalactionex}
  Let $S = \Q[t_1,t_2]$, $R_n = \Q[x_1,\ldots,x_n]$, and consider the
  natural action of ${\mathfrak S}_n$ on $R_n$. Let $Q$ be the kernel
  of the homomorphism induced by the map $\phi\colon R_3 \to S$ given
  by $\phi(x_1) = t_1^2$, $\phi(x_2) = t_2^2$, and $\phi(x_3) = t_1
  t_2$.  Then, $Q = \langle x_1x_2-x_3^2 \rangle$, and $L_4(Q)
  \subseteq R_4$ is generated by the following $12$ polynomials:
\[
\begin{split}
  & x_1x_2-x_3^2,\ x_1x_2-x_4^2,\ x_1x_3-x_2^2,\ x_1x_3-x_4^2, \\
  & x_1x_4-x_3^2,\ x_1x_4-x_2^2,\ x_2x_3-x_1^2,\ x_2x_3-x_4^2, \\
  & x_2x_4-x_1^2,\ x_2x_4-x_3^2,\ x_3x_4-x_1^2,\ x_3x_4-x_2^2. \\
\end{split}
\]
\end{ex}

We would also like a way to project a set of elements in $R_m$ down to
a smaller ring $R_n$ ($n \leq m$).

\begin{defn}
  Let $B \subseteq R_m$ and $n\leq m$.  The \textit{$n$-projection}
  $P_n(B)$ of $B$
  is the $\frak S_n$-invariant ideal of $R_n$ given by \[P_n(B) = \< g
  : g \in B \>_{R_m[{\mathfrak S}_m]} \cap R_n.\]
\end{defn}
We now consider increasing chains $I_\circ$ of ideals $I_n \subseteq
R_n$: \[ I_1 \subseteq I_2 \subseteq \cdots \subseteq I_n \subseteq
\cdots,\] simply called {\em chains}\/ below. Of course, such chains
will usually fail to stabilize since they are ideals in larger and
larger rings.  However, it is possible for these ideals to stabilize
``up to the action of the symmetric group'', a concept we make clear
below.  For the purposes of this work, we will only consider a special
class of chains; namely, a \textit{symmetrization invariant chain}
(resp.~\textit{projection invariant chain}) is one for which $L_{m}(I_n)
\subseteq I_{m}$ (resp.~$P_{n}(I_m) \subseteq I_{n}$) for all $n\leq m$.  If
$I_\circ$ is both a symmetrization and a projection invariant chain,
then it will be simply called an \textit{invariant chain}.  We will
encounter some concrete invariant chains in the next section.  The
stabilization definition alluded to above is as follows.

\begin{defn}\label{stabdef}
  A symmetrization invariant chain of ideals $I_\circ$ as above
  \textit{stabilizes modulo the symmetric group} (or simply
  \textit{stabilizes}) if there exists a positive integer $N$ such
  that \[L_{m}(I_n) = I_m \qquad\text{for all $m \geq n > N$.}\]
\end{defn}

To put it another way, accounting for the natural action of the
symmetric group, the ideals $I_n$ are the same for large enough $n$.
Let us remark that if for a symmetrization invariant chain $I_\circ$,
there is some integer $N$ such that $L_m(I_N) = I_m$ for all $m > N$,
then $I_\circ$ stabilizes. This follows from the inclusions \[I_m =
L_m(I_N) \subseteq L_m(I_n)\subseteq I_m, \ \ n > N.\] Any chain
$I_\circ$ naturally gives rise to an ideal $\mathcal{I}(I_\circ)$ of
$R = \bigcup_{n\geq 1} R_n$ by way of \[ \mathcal{I}(I_\circ) :=
\bigcup_{n \geq 1} I_n.\] Conversely, if $I$ is an ideal of $R$, then
\[I_n = \mathcal{J}_n(I) := I \cap R_n\] defines the components of a
chain $\mathcal{J}(I) := I_\circ$.  Clearly, for any ideal $I
\subseteq R$, we have $\mathcal{I} \circ \mathcal{J}(I) = I$, but, as
is easily seen, it is not true in general that $\mathcal{J} \circ
\mathcal{I}(I_\circ) = I_\circ$.  However, for invariant chains, this
relationship does hold, as the following straightforward lemma
describes.

\begin{lem}
  There is a one-to-one, inclusion-preserving correspondence between
  invariant chains $I_\circ$ and invariant ideals $I$ of $R$ given by
  the maps $\mathcal{I}$ and $\mathcal{J}$. \qed
\end{lem}

For the remainder of this section we consider the case where, for a
commutative Noetherian ring $A$, we have $R_n = A[x_1,\ldots,x_n]$ for
each $n$, endowed with the natural action of ${\mathfrak S}_n$ on the
indeterminates $x_1,\dots,x_n$.  Then $R=A[X^\diamond]$ where
$X=\{x_1,x_2,\dots\}$.  We use the results of the previous section to
demonstrate the following.

\begin{thm}\label{symminvarstab}
  Every symmetrization invariant chain stabilizes modulo the symmetric
  group.
\end{thm}

\begin{proof}
  Given a symmetrization invariant chain, construct the invariant
  ideal $I = \mathcal{I}(I_\circ)$ of $R$.  One would now like to
  apply Theorem~\ref{onevarfinitegenthm}; however, more care is needed
  to prove stabilization.  Let $\leq$ be a well-ordering of $X$ of
  order type $\omega$, and let $B$ be a finite Gr\"obner basis for $I$
  with respect to the corresponding term ordering $\leq_{\lex}$ of
  $X^\diamond$ (Theorem~\ref{wellquasithm2} and
  Proposition~\ref{finite GB}). Choose a positive integer $N$ such
  that $B\subseteq I_N$; we claim that $I_m=L_m(I_N)$ for all $m\geq
  N$. Let $f\in I_m$, $f\neq 0$. By the equivalence of (1) and (3) in
  Lemma~\ref{char GB} we have
  $f\underset{B}{\overset{*}\longrightarrow} 0$.  Hence by
  Lemma~\ref{reduction} there are $g_1,\dots,g_n\in B$,
  $h_1,\dots,h_n\in R$, as well as $\sigma_1,\dots,\sigma_n\in \frak
  S_\infty$, such that $$f=h_1\sigma_1 g_1+\cdots+h_n\sigma_n g_n
  \quad\text{and}\quad \lm(f)=\max_i \lm(h_i\sigma_i g_i).$$
  By
  Remark~\ref{reduction-remark, 2} we may assume that in fact
  $\sigma_i\in\frak S_m$ for each $i$.  Moreover $\lm(h_i)\leq_{\lex}
  \lm(f)$; hence $|\lm(h_i)|\leq |\lm(f)|\leq m$, for each $i$.
  Therefore $h_i\in R_m$ for each $i$.  This shows that $f\in
  L_m(B)\subseteq L_m(I_N)$ as desired.
\end{proof}

\section{A Chemistry Motivation}\label{chemmotivation}

We can now discuss the details of the basic problem that is of
interest to us.  It was brought to our attention by Bernd Sturmfels,
who, in turn, learned about it from Andreas Dress.

Fix a natural number $k\geq 1$. Given a set $S$ we denote by $\<S\>^k$
the set of all ordered $k$-element subsets of $S$; that is, $\<S\>^k$
is the set of all $k$-tuples $\u=(u_1,\dots,u_k)\in S^k$ with pairwise
distinct $u_1,\dots,u_k$. We also just write $\<n\>^k$ instead of
$\<\{1,\dots,n\}\>^k$.  Let $K$ be a field, and for $n\geq k$ consider
the polynomial ring $$R_n = K\big[\{x_{\u}\}_{\u\in\<n\>^k}\big].$$
We
let $\frak S_n$ act on $\<n\>^k$ by
$$\sigma (u_1, \ldots, u_k) = \big(\sigma(u_1),\ldots,
\sigma(u_k)\big).$$
This induces an action $(\sigma,x_{\u})\mapsto
\sigma x_{\u}=x_{\sigma\u}$ of $\frak S_n$ on the indeterminates
$x_\u$, which we extend to an action of $\frak S_n$ on $R_n$ in the
natural way.  We also put $R=\bigcup_{n\geq k} R_n$. Note that $$R=
K\big[\{x_{\u}\}_{\u\in\<\Omega\>^k}\big],$$
where
$\Omega=\{1,2,3,\dots\}$ is the set of positive integers, and that the
actions of $\frak S_n$ on $R_n$ combine uniquely to an action of
$\frak S_\infty$ on $R$.  Now let $f(y_1,\ldots,y_k) \in
K[y_1,\ldots,y_k]$, let $t_1,t_2,\dots$ be an infinite sequence of
pairwise distinct indeterminates over $K$, and for $n \geq k$ consider
the $K$-algebra homomorphism
\[
\phi_n\colon R_n \to K[t_1,\ldots,t_n], \qquad x_{(u_1, \ldots, u_k)}
\mapsto f(t_{u_1},\ldots,t_{u_k}).
\]
The ideal \[Q_n = \text{ker} \ \phi_n\] of $R_n$ determined by such a
map is the prime ideal of algebraic relations between the quantities
$f(t_{u_1},\ldots,t_{u_k})$.  Such ideals arise in chemistry
\cite{Ruch1, Ruch2, Ruch3}; of specific interest is when $f$ is
a Vandermonde polynomial $\prod_{i<j} (y_i-y_j)$.  In this case, the
ideals $Q_n$ correspond to relations among a series of experimental
measurements.  One would then like to understand the limiting behavior
of such relations, and in particular, to see that they stabilize up to
the action of the symmetric group.

\begin{ex}
  The permutation $\sigma = (1\,2\,3) \in {\mathfrak S}_3$ acts on the
  elements
\[(1,2),\ (2,1),\ (1,3),\ (3,1),\ (2,3),\ (3,2)\] of $\<3\>^2$ to give
\[(2,3),\ (3,2),\ (2,1),\ (1,2),\ (3,1),\ (1,3),\]
respectively.  Let $f(t_1,t_2) = t_1^2t_2$.  Then the action of
$\sigma$ on the valid relation $x_{12}^2x_{31}-x_{13}^2x_{21} \in Q_3$
gives us another relation $x_{23}^2x_{12}-x_{21}^2x_{32}\in Q_3$.
\end{ex}

It is easy to see that, by construction, the chain $Q_\circ$ of ideals
\[ Q_k \subseteq Q_{k+1} \subseteq \cdots \subseteq Q_n \subseteq
\cdots\] (which we call the chain of ideals {\it induced by the
  polynomial $f$}\/) is an invariant chain.  As in the proof of
Theorem~\ref{symminvarstab}, we would like to form the ideal $Q =
\bigcup_{n \geq k} Q_n$ of the infinite-dimensional polynomial ring $R
= \bigcup_{n \geq k} R_n$, and then apply a finiteness theorem to
conclude that $Q_\circ$ stabilizes in the sense mentioned above
(Definition \ref{stabdef}).  For $k=1$, Theorem~\ref{symminvarstab}
indeed does the job.  Unfortunately however, this simple-minded
approach fails for $k\geq 2$:

\begin{prop}\label{counterexample}
  For $k\geq 2$, the $R[{\mathfrak S}_{\infty}]$-module $R$ is not
  Noetherian.
\end{prop}
\begin{proof}
  Let us make the dependence on $k$ explicit and denote $R$ by
  $R^{(k)}$.  Then $$x_{(u_1,\dots,u_k,u_{k+1})}\mapsto
  x_{(u_1,\dots,u_k)}$$
  defines a surjective $K$-algebra
  homomorphism $\pi_k\colon R^{(k+1)}\to R^{(k)}$ with invariant
  kernel. Hence if $R^{(k+1)}$ is Noetherian as an $R[{\mathfrak
    S}_{\infty}]$-module, then so is $R^{(k)}$; thus it suffices to
  prove the proposition in the case $k=2$.  Suppose therefore that
  $k=2$. By Lemma~\ref{converse of finitenessthm} it is enough to
  produce an infinite bad sequence for the quasi-ordering $|_{\frak S_\infty}$
  of $X^\diamond$, where
  $X=\{x_{\i}:\i\in\<\Omega\>^2\}$. For this, consider the sequence of
  monomials
\begin{align*}
  s_3 &= x_{(1,2)}x_{(3,2)}x_{(3,4)} \\
  s_4 &= x_{(1,2)}x_{(3,2)}x_{(4,3)}x_{(4,5)} \\
  s_5 &= x_{(1,2)}x_{(3,2)}x_{(4,3)}x_{(5,4)}x_{(6,7)} \\
  &\vdots \\
  s_n &= x_{(1,2)}x_{(3,2)}x_{(4,3)}\cdots x_{(n,n-1)}x_{(n,n+1)}
  \qquad (n=3,4,\dots) \\
  &\vdots
\end{align*}
Now for $n<m$ and any $\sigma\in\frak S_\infty$, the monomial $\sigma
s_n$ does not divide $s_m$. To see this, suppose otherwise. Note that $x_{(1,2)}$,
$x_{(3,2)}$ is the only pair of indeterminates which divides $s_n$ or
$s_m$ and has the form $x_{(i,j)}$, $x_{(l,j)}$ ($i,j,l\in\Omega$).
Therefore $\sigma(2)=2$, and either $\sigma(1)=1$, $\sigma(3)=3$, or
$\sigma(1)=3$, $\sigma(3)=1$. But since $1$ does not appear as the
second component $j$ of a factor $x_{(i,j)}$ of $s_m$, we have
$\sigma(1)=1$, $\sigma(3)=3$. Since $x_{(4,3)}$ is the only
indeterminate dividing $s_n$ or $s_m$ of the form $x_{(i,3)}$ with
$i\in\Omega$, we get $\sigma(4)=4$; since $x_{(5,4)}$ is the only
indeterminate dividing $s_n$ or $s_m$ of the form $x_{(i,4)}$ with
$i\in\Omega$, we get $\sigma(5)=5$; etc. Ultimately this yields
$\sigma(i)=i$ for all $i=1,\dots,n$.  But the only indeterminate
dividing $s_m$ of the form $x_{(n,j)}$ with $j\in\Omega$ is
$x_{(n,n-1)}$; hence the factor $\sigma
x_{(n,n+1)}=x_{(n,\sigma(n+1))}$ of $\sigma s_n$ does not divide
$s_m$.  This shows that $s_3,s_4,\dots$ is a bad sequence for the
quasi-ordering $|_{\frak S_\infty}$, as claimed.
\end{proof}

\begin{rem}
  The construction of the infinite bad sequence $s_3,s_4,\dots$ in the
  proof of the previous proposition was inspired by an example in
  \cite{Jenkyns-NW}.
\end{rem}

\subsection{A criterion for stabilization}
Our next goal is to give a condition for the chain $Q_\circ$ to
stabilize.  Given $g \in R$, we define the \textit{variable size} of
$g$ to be the number of distinct indeterminates $x_\u$ that appear in
$g$.  For example, $g = x_{12}^5+x_{45}x_{23}+x_{45}$ has variable
size $3$.

\begin{lem}\label{sizesuffcond}
  A chain of ideals $Q_\circ$ induced by a polynomial $f\in
  K[y_1,\dots,y_k]$ stabilizes modulo the symmetric group if and only
  if there exist integers $M$ and $N$ such that for all $n > N$, there
  are generators for $Q_n$ with variable sizes at most $M$. Moreover,
  in this case a bound for stabilization is given by
  ${\rm{max}}(N,kM)$.
\end{lem}

\begin{proof}
  Suppose $M$ and $N$ are integers with the stated property. To see
  that $Q_\circ$ stabilizes, since $Q_\circ$ is an invariant chain, we
  need only verify that $N' = \text{max}(N,kM)$ is such that $Q_m
  \subseteq L_m(Q_n)$ for $m \geq n > N'$.  For this inclusion, it
  suffices that each generator in a generating set for the ideal $Q_m$
  of $R_m$ is in $L_m(Q_n)$.  Since $m > N$, there are generators $B$
  for $Q_m$ with variable sizes at most $M$.  If $g \in B$, then there
  are at most $kM$ different integers appearing as subscripts of
  indeterminates in $g$.  We can form a permutation $\sigma \in {\mathfrak
    S}_m$ such that $\sigma g \in R_{N'}$ and thus in $R_n$.  But then
  $\sigma g \in P_n(Q_m) \subseteq Q_n$ so that $g = \sigma^{-1}\sigma
  g \in L_m(Q_n)$ as desired.

  Conversely, suppose that $Q_{\circ}$ stabilizes.  Then there exists
  an $N$ such that $Q_m = L_m(Q_N)$ for all $m > N$. Let $B$ be any
  finite generating set for $Q_N$.  Then for all $m > N$, $Q_m =
  L_m(B)$ is generated by elements of bounded variable size by
  Lemma~\ref{symmetrgeneratorslemma}.
\end{proof}

Although this condition is a very simple one, it will prove useful.
Below we will apply it together with a preliminary reduction to the
case that each indeterminate $y_1,\dots,y_k$ actually occurs in the
polynomial $f$, which we explain next.  For this we let $\pi_k\colon
R^{(k+1)}\to R^{(k)}$ be the surjective $K$-algebra homomorphism
defined in the proof of Proposition~\ref{counterexample}.  We write
$Q^{(k)}$ for $Q$, and considering $f\in K[y_1,\dots,y_k]$ as an
element of $K[y_1,\dots,y_k,y_{k+1}]$, we also let $Q^{(k+1)}$ be the
kernel of the $K$-algebra homomorphism \begin{multline*} R^{(k+1)} \to
  K[t_1,t_2,\ldots], \qquad x_{(u_1, \ldots, u_k,u_{k+1})} \mapsto
  f(t_{u_1},\ldots,t_{u_k},t_{u_{k+1}}) \\
  (=f(t_{u_1},\ldots,t_{u_k})).
\end{multline*}
Note that $\pi_k(Q^{(k+1)})=Q^{(k)}$, and the ideal $\ker\pi_k$ of
$R^{(k+1)}$ is generated by the elements
$$x_{(u_1,\dots,u_k,i)} - x_{(u_1,\dots,u_k,j)} \qquad
(i,j\in\Omega);$$
in particular, $\ker\pi_k\subseteq Q^{(k+1)}$.  It is
easy to see that as an $R^{(k+1)}[\frak S_\infty]$-module, $\ker\pi_k$ is
generated by the single element
$x_{(1,\dots,k,k+1)}-x_{(1,\dots,k,k+2)}$.  These observations now
yield:

\begin{lem}
  Suppose that the invariant ideal $Q^{(k)}$ of $R^{(k)}$ is finitely
  generated as an $R^{(k)}[\frak S_\infty]$-module. Then the invariant
  ideal $Q^{(k+1)}$ of $R^{(k+1)}$ is finitely generated as
  an $R^{(k+1)}[\frak S_\infty]$-module. \qed
\end{lem}

We let $\frak S_k$ act on $\<\Omega\>^k$ by $$\tau
(u_1,\dots,u_k)=(u_{\tau(1)},\dots,u_{\tau(k)})\qquad\text{ for
  $\tau\in\frak S_k$, $(u_1,\dots,u_k)\in\<\Omega\>^k$.}$$
This action
gives rise to an action of $\frak S_k$ on
$\{x_{\u}\}_{\u\in\<\Omega\>^k}$ by $\tau x_{\u} = x_{\tau\u}$, which
we extend to an action of $\frak S_k$ on $R$ in the natural way. We
also let $\frak S_k$ act on $K[y_1,\dots,y_k]$ by $\tau
f(y_1,\dots,y_k) = f(y_{\tau(1)},\dots,y_{\tau(k)})$.  Note that
$$\tau Q_k \subseteq \tau Q_{k+1} \subseteq \cdots \subseteq \tau Q_n
\subseteq\cdots$$
is the chain induced by $\tau f$. Using the lemma
above we obtain:

\begin{cor}\label{every variable occurs}
  Let $f\in K[y_1,\dots,y_k]$. There are $i\in\{0,\dots,k\}$ and
  $\tau\in\frak S_k$ such that $\tau f\in K[y_1,\dots,y_{i}]$ and each
  of the indeterminates $y_1,\dots,y_{i}$ occurs in $\tau f$.  If the
  chain of ideals induced by the polynomial $\tau f$ stabilizes, then
  so does the chain of ideals induced by $f$. \qed
\end{cor}

\subsection{Chains induced by monomials}
If the given polynomial $f$ is a monomial, then the homomorphism
$\phi_n$ from above produces a (homogeneous) toric kernel $Q_n$. In
particular, there is a finite set of binomials that generate $Q_n$
(see \cite{ConvPoly}).  Although a proof for the general toric case
eludes us, we do have the following.

\begin{thm}\label{toricsqfreethm}
  The sequence of kernels induced by a square-free monomial $f \in
  K[y_1,\ldots,y_k]$ stabilizes modulo the symmetric group.  Moreover,
  a bound for when stabilization occurs is $N = 4k$.
\end{thm}

To prepare for the proof of this result, we discuss in detail the
toric encoding associated to our problem (see \cite[Chapter
14]{ConvPoly} for more details).  By Corollary~\ref{every variable
  occurs}, we may assume that $f=y_1\cdots y_k$.  Then $g-\tau g\in Q$
for all $g\in R$.  We say that $\u=(u_1,\dots,u_k)\in\<\Omega\>^k$ is
\textit{sorted} if $u_1<\cdots< u_k$, and \textit{unsorted} otherwise;
similarly we say that $x_{\u}$ is sorted (unsorted) if $\u$ is sorted
(unsorted, respectively).  For example, $x_{135}$ is a sorted
indeterminate, whereas $x_{315}$ is not. Consider the set of vectors
\[
\mathcal{A}_n = \big\{(i_1,\ldots,i_n) \in \mathbb Z^n : i_1+\cdots +
i_n = k, \ 0 \leq i_1,\ldots,i_n \leq 1\big\}.
\]
View $\mathcal{A}_n$ as an $n$-by-$\binom{n}{k}$ matrix with entries
$0$ and $1$, whose columns are indexed by sorted indeterminates
$x_{\u}$ and whose rows are indexed by $t_i$ ($i = 1,\ldots,n$). (See
Example \ref{sqfreeex} below.)  Let $\sort(\,\cdot\,)$ denote the
operator which takes any word in $\{1,\dots,n\}^*$ and sorts it in
increasing order.  By \cite[Remark 14.1]{ConvPoly}, the toric ideal
$I_{\mathcal{A}_n}$ associated to $\mathcal{A}_n$ is generated (as a 
$K$-vector space) by the binomials $x_{\u_1}\cdots x_{\u_r} -
x_{\v_1}\cdots x_{\v_r}$, where $r\in\N$ and the $\u_i$, $\v_j$ are
sorted elements of $\<n\>^k$ such that
$\sort(\u_1\cdots\u_r)=\sort(\v_1\cdots\v_r)$. In particular, we have
$I_{\cal A_n}\subseteq Q_n$.  Let $B$ be any set of generators for the
ideal $I_{\cal A_n}$.

\begin{lem}\label{sqfreetoriclemma}
  A generating set for the ideal $Q_n$ of $R_n$ is given by \[S = B
  \cup \{x_{\u} - x_{\tau \u} : \ \tau \in {\mathfrak S}_k, \ \u
  \text{ is sorted} \}.\]
\end{lem}

\begin{proof}
  Elements of $Q_n$ are of the form $g = x_{\u_1} \cdots x_{\u_r} -
  x_{\v_1} \cdots x_{\v_r}$, in which the ${\u_i}$ and ${\v_j}$ are
  ordered $k$-element subsets of $\{1,\ldots,n\}$ such that
  $\sort(\u_1\cdots\u_r)=\sort(\v_1\cdots\v_r)$.  We induct on the
  number $t$ of ${\u_i}$ and ${\v_j}$ that are not sorted.  If $t =
  0$, then $g \in I_{\mathcal{A}_n}$, and we are done.  Suppose now
  that $t > 0$ and assume without loss of generality that ${\u_1}$ is
  not sorted.  Let $\tau \in {\mathfrak S}_k$ be such that $\tau
  {\u_1}$ is sorted, and consider the element
  $h=x_{\tau\u_1}x_{\u_2}\cdots x_{\u_r} - x_{\v_1} \cdots x_{\v_r}$
  of $Q_n$. This binomial involves $t - 1$ unsorted indeterminates,
  and therefore, inductively, can be expressed in terms of $S$. But
  then
  \[g = h -
  (x_{\tau \u_1} - x_{\u_1})x_{\u_2} \cdots x_{\u_r}\] can as well,
  completing the proof.
\end{proof}

\begin{ex}\label{sqfreeex}
  Let $k = 2$ and $n = 4$.  Then
\[\begin{array}{*{20}c}
   {} & {x_{12} } & {x_{13} } & {x_{14} } & {x_{23} } & {x_{24} } & {x_{34} }  \\
   {t_1 } & 1 & 1 & 1 & 0 & 0 & 0  \\
   {t_2 } & 1 & 0 & 0 & 1 & 1 & 0  \\
   {t_3 } & 0 & 1 & 0 & 1 & 0 & 1  \\
   {t_4 } & 0 & 0 & 1 & 0 & 1 & 1  \\
 \end{array}\] represents the matrix associated
 to $\mathcal{A}_4$. The ideal $I_{\mathcal{A}_4}$ is generated by the
 two binomials $x_{13}x_{24}-x_{12}x_{34}$ and
 $x_{14}x_{23}-x_{12}x_{34}$.  Hence $Q_4$ is generated by these two
 elements along with
\[\{x_{12}-x_{21},x_{13}-x_{31},x_{14}-x_{41},x_{23}-x_{32},x_{24}-x_{42},x_{34}-x_{43}
\}.\]
\end{ex}

We are now in a position to prove Theorem~\ref{toricsqfreethm}.

\begin{proof}[Proof of Theorem \ref{toricsqfreethm}]
  By Lemma \ref{sizesuffcond}, we need only show that there exist
  generators for $Q_n$ which have bounded variable sizes.  Using
  \cite[Theorem 14.2]{ConvPoly}, it follows that $I_{\mathcal{A}_n}$
  has a quadratic (binomial) Gr\"{o}bner basis for each $n$ (with
  respect to some term ordering of $R_n$).  By Lemma
  \ref{sqfreetoriclemma}, there is a set of generators for $Q_n$ with
  variable sizes at most $4$.  This proves the theorem.
\end{proof}

We close with a conjecture that generalizes Theorem
\ref{toricsqfreethm}.

\begin{conj}\label{toriconj}
  The sequence of kernels induced by a monomial $f$ stabilizes modulo
  the symmetric group.
\end{conj}

\section{Acknowledgment}
We would like to thank Bernd Sturmfels for bringing the problem found
in Section \ref{chemmotivation} (originating from Andreas Dress) to
our attention and for making us aware of Theorem 14.2 in
\cite{ConvPoly}.

\end{document}